\documentclass[12pt]{article} 
\usepackage{amssymb,latexsym,color,psfig}
\newcounter{environment}[section]
\renewcommand{\theenvironment}{%
\arabic{section}.\arabic{environment}}

{\begin{rm}\refstepcounter{environment}{\textbf\theenvironment\
\bf Definition.~~}}%
{\end{rm}}

{\begin{rm}\refstepcounter{environment}{\textbf\theenvironment\
\bf Example.~~}}%
{\end{rm}}

{\begin{rm}\refstepcounter{environment}{\textbf\theenvironment\
\bf Conjecture.~~}}%
{\end{rm}}

\newenvironment{proposition}%
{\begin{rm}\refstepcounter{environment}{\textbf\theenvironment\
\bf Proposition.~~}}%
{\end{rm}}

{\begin{rm}\refstepcounter{environment}{\textbf\theenvironment\
\bf Proposition and Definition.~~}}%
{\end{rm}}

\newenvironment{theorem}%
{\begin{rm}\refstepcounter{environment}{\textbf\theenvironment\
\bf Theorem.~~}}%
{\end{rm}}

\newenvironment{corollary}%
{\begin{rm}\refstepcounter{environment}{\textbf\theenvironment\
\bf Corollary.~~}}%
{\end{rm}}

\newenvironment{lemma}%
{\begin{rm}\refstepcounter{environment}{\textbf\theenvironment\
\bf Lemma.~~}}%
{\end{rm}}

\begin{document}
\newcommand{\qbc}[2]{ {\left [{#1 \atop #2}\right ]}}
\newcommand{\anbc}[2]{{\left\langle {#1 \atop #2} \right\rangle}}
\newcommand{\be}{\begin{enumerate}}
\newcommand{\ee}{\end{enumerate}}
\newcommand{\beq}{\begin{equation}}
\newcommand{\eeq}{\end{equation}}
\newcommand{\bea}{\begin{eqnarray}}
\newcommand{\eea}{\end{eqnarray}}
\newcommand{\beas}{\begin{eqnarray*}}
\newcommand{\eeas}{\end{eqnarray*}}
\newcommand{\zz}{\mathbb{Z}}
\newcommand{\pp}{\mathbb{P}}
\newcommand{\nn}{\mathbb{N}}
\newcommand{\qq}{\mathbb{Q}}
\newcommand{\rr}{\mathbb{R}}
\newcommand{\bm}[1]{{\mbox{\boldmath $#1$}}}
\newcommand{\bms}[1]{{\mbox{\scriptsize\boldmath $#1$}}}
\newcommand{\bmt}{\bm{T}}
\newcommand{\bmd}{\bm{D}}
\newcommand{\rk}{\mathrm{rank}}
\newcommand{\sse}{\mathrm{SS}}
\newcommand{\co}{\mathrm{code}}
\newcommand{\rco}{\overline{\mathrm{code}}}
\newcommand{\hgt}{\mathrm{ht}}
\newcommand{\htt}{\mathrm{ht}(\bm{T})}
\newcommand{\lm}{\lambda/\mu}
\newcommand{\rc}{\mathrm{rc}}
\newcommand{\fo}{\mathfrak{o}}
\newcommand{\sn}{\mathfrak{S}_n}
\newcommand{\st}{\,:\,}
\newcommand{\tr}{\textcolor{red}}
\newcommand{\tb}{\textcolor{blue}}
\newcommand{\tg}{\textcolor{green}}
\newcommand{\tm}{\textcolor{magenta}}
\newcommand{\tbn}{\textcolor{brown}}
\newcommand{\tp}{\textcolor{purple}}
\newcommand{\tn}{\textcolor{nice}}
\newcommand{\tor}{\textcolor{orange}}

\definecolor{brown}{cmyk}{0,0,.35,.65}
\definecolor{purple}{rgb}{.5,0,.5}
\definecolor{nice}{cmyk}{0,.5,.5,0}
\definecolor{orange}{cmyk}{0,.35,.65,0}

\begin{centering}
\textcolor{red}{\Large\bf The Rank and Minimal Border Strip}\\
\textcolor{red}{\Large\bf Decompositions  of a Skew Partition}\\[.2in] 
\textcolor{blue}{Richard P. Stanley}\footnote{Partially supported by
  NSF grant \#DMS-9988459 and by the Isaac Newton Institute for
  Mathematical Sciences.}\\ 
Department of Mathematics\\
Massachusetts Institute of Technology\\
Cambridge, MA 02139\\
\emph{e-mail:} rstan@math.mit.edu\\[.2in]
\textcolor{magenta}{version of 13 September 2001}\\[.2in]

\end{centering}
\vskip 10pt
\section{Introduction.} \label{sec1}
Let $\lambda=(\lambda_1,\lambda_2,\dots)$ be a
partition of the integer $n$, i.e., $\lambda_1\geq \lambda_2\geq
\cdots \geq 0$ and $\sum \lambda_i=n$. The (Durfee or Frobenius)
\emph{rank} of 
$\lambda$, denoted rank$(\lambda)$, is the length of the main diagonal
of the diagram of $\lambda$, or equivalently, the largest integer $i$
for which $\lambda_i\geq i$ \cite[p.\ 289]{ec2}. We
will assume familiarity with the notation and terminology involving
partitions and symmetric functions found in \cite{macd} and
\cite{ec2}. Nazarov and Tarasov \cite[{\S}1]{naz}, in connection with tensor
products of Yangian modules $Y(\mathfrak{gl}_n)$, defined a 
generalization of rank to skew partitions (or skew diagrams)
$\lambda/\mu$. There are several simple equivalent definitions of
rank$(\lm)$ which we summarize in Proposition~\ref{prop1}. In
particular, rank$(\lambda/\mu)$ is the least integer $r$ such that
$\lambda/\mu$ is a disjoint union of $r$ border strips (also called
ribbons or rim hooks). In Section~\ref{sec4} we consider the structure
of the decompositions of $\lambda/\mu$ into this minimal number $r$ of
border strips. For instance, we show that the number of ways to write
$\lm$ as a disjoint union of $r$ border strips is a perfect
square. A consequence of our results will be that if $\chi^{\lm}$ is
the skew character of the symmetric group $\sn$ indexed by $\lm$ and
if $w$ is a permutation in $\sn$ with $\rk(\lm)$ cycles (in its
disjoint cycle decomposition) for which exactly $m_i$ cycles have
length $i$, then $\chi^{\lm}(w)$ is divisible by
$m_1(w)!m_2(w)!\cdots$.

In addition to the various characterizations of rank$(\lambda/\mu)$
given by Proposition~\ref{prop1} we have a further possible
characterization which we have been unable to prove or disprove.
Namely, let $s_{\lambda/\mu}(1^t)$ denote the skew Schur function
$s_{\lambda/\mu}$ evaluated at $x_1=\cdots=x_t=1$, $x_i=0$ for $i>t$.
For fixed $\lambda/\mu$, $s_{\lambda/\mu}(1^t)$ is a polynomial in
$t$. Let zrank$(\lambda/\mu)$ denote the exponent of the largest power
of $t$ dividing $s_{\lambda/\mu}(1^t)$ (as a polynomial in $t$). It is
easy to see (Proposition~\ref{prop:zrank}) that zrank$(\lm)\geq
\rk(\lm)$, and we ask whether equality always holds. We know of two
main cases where the answer is affirmative: (1) when $\lambda/\mu$
is an ordinary partition (i.e., $\mu=\emptyset$), a trivial
consequence of known results on Schur functions
(Theorem~\ref{thm1}(a)), 
and (2) when every row of the Jacobi-Trudi matrix for $\lm$ which
contains an entry equal to 0 also contains an entry equal to 1
(Theorem~\ref{thm1}(b)).

\section{Characterizations of Frobenius rank.} \label{sec2}
\begin{figure}
\centerline{\psfig{figure=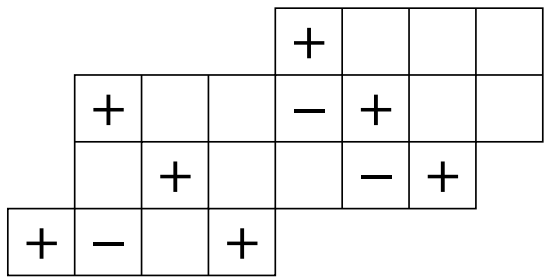}}
\caption{Outside and inside diagonals of the skew shape $8874/411$}
\label{fig1}
\end{figure}

Let $\lm$ be a skew shape, which we identify with its Young
diagram $\{(i,j)\st \mu_i<j\leq \lambda_i\}$. We regard the points
$(i,j)$ of the Young diagram as squares. An \emph{outside top
corner} of $\lm$ is a square $(i,j)\in\lm$ such that $(i-1,j),
(i,j-1)\not\in \lm$. An \emph{outside diagonal} of $\lm$ consists of
all squares $(i+p,j+p)\in\lm$ for which $(i,j)$ is a fixed outside top
corner. Similarly an \emph{inside top corner} of $\lm$ is a square
$(i,j)\in\lm$ such that $(i-1,j),(i,j-1)\in\lm$ but $(i-1,j-1)\not \in
\lm$. An \emph{inside diagonal} of $\lm$ consists of all squares
$(i+p,j+p)\in\lm$ for which $(i,j)$ is a fixed inside top
corner. If $\mu=\emptyset$, then $\lm$ has one outside diagonal (the
main diagonal) and no inside diagonals. Figure~\ref{fig1} shows the
skew shape $8874/411$, with outside diagonal squares marked by + and
inside diagonal squares by $-$. 

Let $d^+(\lm)$ (respectively, $d^-(\lm)$) denote the total number of
outside diagonal squares (respectively, inside diagonal squares) of
$\lm$.  Following Nazarov and Tazarov \cite[{\S}1]{naz}, we define the
(Durfee or Frobenius) \emph{rank} of $\lm$, denoted $\rk(\lm)$, to be
$d^+(\lm)-d^-(\lm)$.  Clearly when $\mu=\emptyset$ this reduces to the
usual definition of $\rk(\lambda)$ mentioned in the introduction. We
see, for instance, from Figure~\ref{fig1} that $\rk(8874/411)=4$.

\begin{figure}
\centerline{\psfig{figure=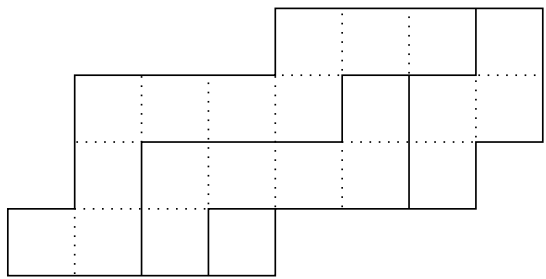}}
\caption{A minimal border strip decomposition of the skew shape
  $8874/411$} 
\label{fig2}
\end{figure}

We wish to give several equivalent definitions of $\rk(\lm)$. First we
discuss the necessary background. A skew shape $\lm$ is
\emph{connected} if the interior of the Young diagram of $\lm$,
regarded as a union of solid squares, is a connected (open) set.  A
\emph{border strip} \cite[p.\ 345]{ec2} is a connected skew shape with
no $2\times 2$ square. (The empty diagram $\emptyset$ is \emph{not} a
border strip.) A border strip is uniquely determined, up to
translation, by its row lengths; there are exactly $2^{n-1}$ border
strips with $n$ squares (up to translation). We say that a border
strip $B\subseteq \lm$ is a \emph{border strip of $\lm$} if $\lm-B$ is
a skew shape $\nu/\mu$ (so $B=\lambda/\nu$). Equivalently, we say that
$B$ can be \emph{removed} from $\lm$. A border strip $B$ of $\lm$ is
determined by its lower left-hand square init$(B)$ and upper
right-hand square fin$(B)$. A \emph{border strip decomposition}
\cite[p.\ 470]{ec2} of $\lm$ is a partitioning of the squares of $\lm$
into (pairwise disjoint) border strips. Let $N=|\lm|:=\sum \lambda_i
-\sum\mu_i$ and $\sigma=
(\sigma_1,\dots,\sigma_\ell)\vdash N$, where $\sigma_\ell>0$.
We say that a border strip decomposition $\bm{D}$ has \emph{type}
$\sigma\vdash N$ if the sizes (number of squares) of the border
strips appearing in $\bm{D}$ are $\sigma_1,\dots,\sigma_\ell$.
A border strip decomposition
of $\lm$ is \emph{minimal} if the number of border strips is
minimized, i.e., there does not exist a border strip decomposition
with fewer border strips.  Figure~\ref{fig2} shows a minimal border
strip decomposition of the skew shape $8874/411$.

A concept closely related to border strip decompositions is that of
border strip tableaux \cite[p.\ 346]{ec2}. Let $\lm\vdash N$. Let
$\alpha=(\alpha_1,\alpha_2,\dots,\alpha_m)$ be a composition of
$N$, i.e., $\alpha_i\in\pp=\{1,2,\dots\}$ and $\sum \alpha_i=N$. A
\emph{border strip tableau} of (shape) $\lm$ and type $\alpha$ is a
sequence 
  \beq \mu=\lambda^0\subset \lambda^1\subset \cdots \subset
    \lambda^r =\lambda \label{eq:bsdef} \eeq
such that $\lambda^i/\lambda^{i-1}$ is a border strip of size
$\alpha_i$. (Note that the 
type of a border strip decomposition is a partition but of a border
strip tableau is a composition.) Every border
strip tableau $\bmt$ of shape $\lm$ defines a border strip
decomposition $\bmd$ of $\lm$, viz., the border strips
$\lambda^i/\lambda^{i-1}\neq\emptyset$ of $\bmt$ are just the border
strips of $\bmd$. We say that $\bmd$ \emph{corresponds} to $\bmt$ and
conversely that $\bmt$ \emph{corresponds} to $\bmd$. Of course given
$\bmt$, the corresponding $\bmd$ is unique, but not conversely. If
$\bmt$ corresponds to a minimal border strip decomposition $\bmd$,
then we call $\bmt$ a \emph{minimal} border strip tableau.


Now suppose that $\ell(\lambda)\leq n$, where $\ell(\lambda)$ denotes the
number of (nonzero) parts of $\lambda$. Recall that the Jacobi-Trudi
identity for the skew Schur function $s_{\lm}$ \cite[Thm.\
7.16.1]{ec2} asserts that
  $$ s_{\lm} = \det\left( h_{\lambda_i-\mu_j-i+j}\right)_{i,j=1}^n, $$
where $h_k$ denotes the complete homogeneous symmetric function of
degree $k$, with the convention $h_0=1$ and $h_k=0$ for $k<0$. Denote
the matrix $\left( h_{\lambda_i-\mu_j-i+j}\right)$ appearing in the
Jacobi-Trudi identity by JT$_{\lm}$, called the \emph{Jacobi-Trudi
matrix} of the skew shape $\lm$. Let jrank$(\lm)$ denote the number of
rows of JT$_{\lm}$ that don't contain a 1. Note that JT$_{\lm}$
implicitly depends on $n$, but jrank$(\lm)$ does not depend on the
choice of $n$.

\begin{figure}
\centerline{\psfig{figure=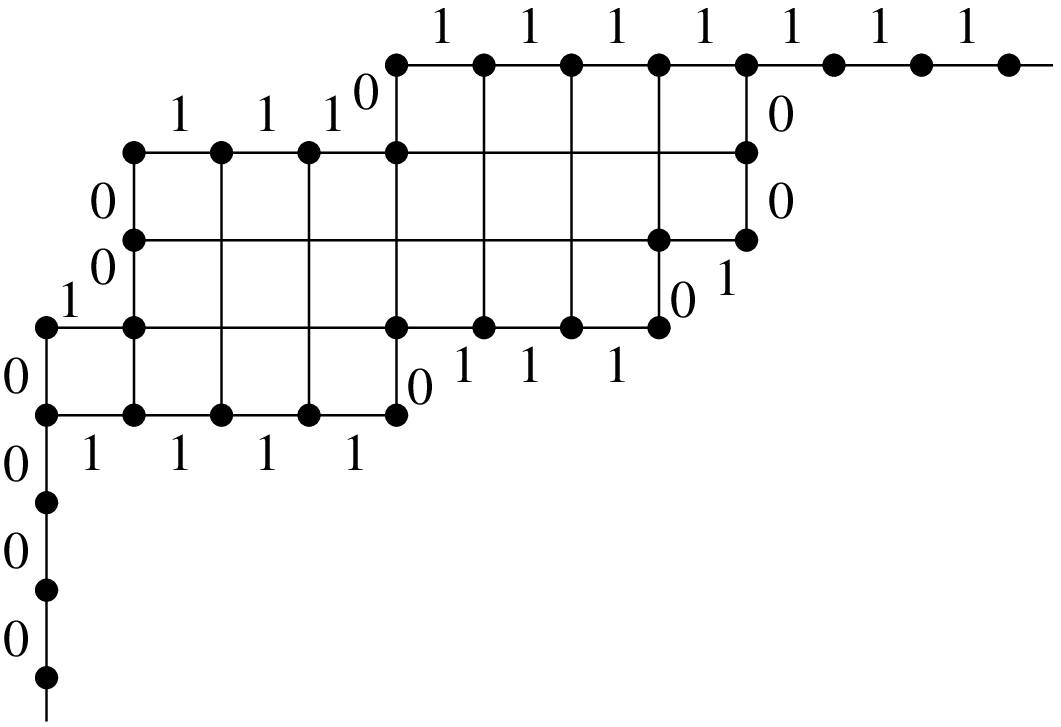}}
\caption{Constructing the code of $8874/411$}
\label{fig3}
\end{figure}

Our final piece of background material concerns the (Com\'et)
code of a shape $\lambda$ \cite[Exer.\ 7.59]{ec2}, generalized
to skew shapes $\lm$.  Let $\lm$ be a skew 
shape, with its left-hand edge and upper edge extended to infinity, as
shown in Figure~\ref{fig3} for $\lm=8874/411$. Put a 0 next to each
vertical edge and a 1 next to each horizontal edge of the ``lower
envelope'' and ``upper envelope'' of $\lm$ (whose definition should be
clear from Figure~\ref{fig3}). If we read these numbers as we move
north and east along the lower envelope we obtain an infinite binary
sequence $C_{\lm}=\cdots c_{-2}c_{-1}c_0c_1c_2\cdots$ beginning with
infinitely many 0's and ending with infinitely many 1's. Similarly if
read these numbers as we move north and east along the upper envelope
we obtain another such binary sequence $D_{\lm}=\cdots d_{-2}d_{-1}d_0
d_1d_2\cdots$. The indexing of the terms of $C_{\lm}$ and $D_{\lm}$ is
arbitrary (it doesn't affect the sequences themselves), but we require
them to ``line up'' in the sense that common steps in the two
envelopes should have common indices. We call the resulting two-line
array
   \beq \co(\lm) = \begin{array}{cccccccc} \cdots & c_{-2} & c_{-1}
  & c_0 & c_1 & c_2 & \cdots\\ \cdots & d_{-2} & d_{-1} & d_0
  & d_1 & d_2 & \cdots \end{array}, \label{eq:comet} \eeq
the (Com\'et) \emph{code} of $\lm$ (also known as the
\emph{partition sequence} of $\lambda$ \cite{bes1}\cite{bes2}). If we
omit the infinitely many initial columns ${0\atop 0}$ and final
columns ${1\atop 1}$ from $\co(\lm)$, then we call the resulting array
the \emph{reduced code} of $\lm$, denoted $\rco(\lm)$. 
Thus for instance from Figure~\ref{fig3} we see that 
  $$ \rco(8874/411) = \begin{array}{cccccccccccc}
   1 & 1 & 1 & 1 & 0 & 1 & 1 & 1 & 0 & 1 & 0 & 0\\
   0 & 1 & 0 & 0 & 1 & 1 & 1 & 0 & 1 & 1 & 1 & 1 \end{array}. $$
\indent A two-line array (\ref{eq:comet}) with infinitely many initial
columns ${0\atop 0}$ and final columns ${1\atop 1}$ is the code of
some $\lm$ if and only if for all $i$,
   \beq \#\{ j\leq i\st (c_j,d_j)=(1,0)\} \geq \#\{ j\leq i\st 
       (c_j,d_j)=(0,1)\}, \label{eq:lp} \eeq
and if 
 \beq \#\{ j\in\mathbb{Z}\st (c_j,d_j)=(1,0)\} = \#\{ j\in\mathbb{Z} 
      \st (c_j,d_j)=(0,1)\}. \label{eq:lp2} \eeq
If $\mu=\emptyset$ then the second row of code$(\lm)$ is redundant, so
we define code$(\lambda)$ to be the first row of code$(\lm)$.
If $\co(\lm)$ is given by (\ref{eq:comet}) then we write $s(c_i)$
(respectively, $s(d_i)$) for the (unique) square of $\lm$ that
contains the edge of the lower envelope (respectively, upper envelope)
of $\lm$ corresponding to $c_i$ (respectively, $d_i$).
The following fundamental property of code$(\lm)$ appears e.g.\ in
\cite[Exer.\ 7.59(b)]{ec2} for ordinary shapes and carries over
directly to skew shapes.

\begin{proposition} \label{prop:codebs}
\emph{Let $\co(\lm)$ be given by (\ref{eq:comet}). Then
removing a border strip of size $p$ from $\lm$ is equivalent to
choosing $i$ with $c_i=1$ and $c_{i+p}=0$, and then replacing $c_i$
with 0 and $c_{i+p}$ with 1, provided that (\ref{eq:lp}) continues to
hold. Specifically, such a pair $(i,i+p)$ corresponds to the border
strip $B$ of size $p$ defined by}
  $$ \mathrm{init}(B)=s(c_i),\quad \mathrm{fin}(B)=s(c_{i+p}). $$
\emph{Moreover,} code$(\lm-B)$ \emph{is obtained from} code$(\lm)$
\emph{by setting $c_i=0$ and $c_{i+p}=1$.} 
\end{proposition}

We can now state several characterizations of $\rk(\lm)$.

\begin{proposition} \label{prop1}
\emph{For any skew shape $\lm$, the following numbers are equal.}
\vspace{-.2in}
\be\item[(a)] rank$(\lm)$
 \item[(b)] \emph{the number of border strips in a minimal border strip
   decomposition of $\lm$}
  \item[(c)] jrank$(\lm)$
  \item[(d)] \emph{the number of columns of $\co(\lm)$
      equal to ${0\atop 1}$ (or to ${1\atop 0}$)}
\ee
\end{proposition}

\textbf{Proof.} 
By equations (\ref{eq:lp}) and (\ref{eq:lp2}) there exists a bijection 
  $$ \vartheta:\{ i\st
  (c_i,d_i)=(1,0)\}\rightarrow \{i\st 
  (c_i,d_i)=(0,1)\} $$ 
such that $\vartheta(i)>i$ for all $i$ in the domain of $\vartheta$.  
By Proposition~\ref{prop:codebs}, as we successively remove border
strips from $\lm$ the bottom line $\cdots d_{-1} d_0 d_1\cdots$ of
code$(\lm)$ remains the same, while the top line $\cdots c_{-1} c_0
c_1\cdots$ interchanges a 0 and 1. We will exhaust all of $\lm$ when
the top line becomes equal to the bottom. Hence the number of
border strips appearing in a border strip decomposition of $\lm$ is at
least the number of columns ${0\atop 1}$ of code$(\lm)$. On the other
hand, we can achieve exactly this number by interchanging $c_i$ with
$c_{\vartheta(i)}$ for all $i$ such that $(c_i,d_i)=(0,1)$. 
It follows that (b) and (d) are equal.

Let $B$ be the (unique) largest border strip of $\lm$ such that
init$(B)$ is the bottom square of the 
leftmost column of $\lm$. $B$ will intersect each diagonal (running from
upper-left to lower-right) of its connected component $\sigma$ of
$\lm$ exactly once. The number of 
outside diagonals of $\sigma$ is one more than the number of inside
diagonals. Hence rank$(\lm)=\rk(\lm-B)+1$. Continuing to remove the
largest border strip results in a minimal border strip decomposition
of $\lm$. (Minimality is an easy consequence of
Proposition~\ref{prop:codebs}.) Since each border strip removal
reduces the rank by one, it follows that (a) and (b) are equal. 

Finally consider the Jacobi-Trudi matrix JT$_{\lm}$. We prove by
induction on the number of rows of JT$_{\lm}$ that (b) and (c) are
equal. The assertion is clear when JT$_{\lm}$ has one row, so assume
that JT$_{\lm}$ has more than one row. We may assume that $\lm$ has no
empty rows, since ``compressing'' $\lm$ by removing all empty rows
does not change (c). Let JT$^\prime_{\lm}$ denote JT$_{\lm}$ with the
first row and last column removed. Let $\nu$ be the shape obtained by
removing a maximal border strip from each connected component of $\lm$
and deleting the bottom (empty) row. If $\lm$ has $c$ connected
components, then $\rk(\nu/\mu)=\rk(\lm)-c$. Now the $(i,j)$-entry
$h_{\lambda_{i+1}-\mu_j-i+j-1}$ of the matrix JT$^\prime_{\lm}$
satisfies
  $$ h_{\lambda_{i+1}-\mu_j-i+j-1} = \left\{ \begin{array}{rl}
    h_{\nu_i-\mu_j-i+j}, & \mbox{if row $i$ of $\lm$ is not the last row 
      of a}\\ & \mbox{ connected component of $\lm$}\\
   h_{\nu_i-\mu_j-i+j+1}, & \mbox{otherwise.} \end{array} \right. $$  
Moreover, if row $i$ is the last row of a connected component of $\lm$
(other than the bottom row of $\lm$) then the $(i,i)$-entry of
JT$_{\nu/\mu}$ is 1, while the $i$th row of JT$_{\lm}$ does not
contain a 1. It follows that jrank$(\nu/\mu)=\mathrm{jrank}(\lm)-c$,
and the equality of (b) and (c) follows by induction. $\ \Box$

The equivalence of (a) and (c) in Proposition~\ref{prop1} is also an
immediate consequence of \cite[Prop.\ 1.32]{naz}.  

The following corollary was first proved by Nazarov and Tarasov
\cite[Thm.\ 1.4]{naz} using the definition
rank$(\lm)=d^+(\lm)-d^-(\lm)$. The result is not obvious (even for
nonskew shapes $\lambda$) using this definition, but it is an
immediate consequence of parts (b) or (d) of Proposition~\ref{prop1}.

\begin{corollary} \label{cor:r}
\emph{Let $(\lm)^\natural$ denote the skew shape obtained by rotating
  the diagram of $\lm$ $180^\circ$, i.e, replacing $(i,j)\in\lm$ with
  $(h-i,k-i)$ for some $h$ and $k$. Then}
  rank$(\lm)=\rk((\lm)^\natural)$.  
\end{corollary}

\section{An open characterization of rank$(\lm)$} \label{sec3}
Recall that in Section~\ref{sec1} we defined zrank$(\lm)$ to be the
largest power of $t$ dividing the polynomial $s_{\lm}(1^t)$.

\textbf{Open problem.} Is it true that 
  \beq \rk(\lm)=\mathrm{zrank}(\lm) \label{eq:mainop} \eeq
for all $\lm$?

\begin{proposition} \label{prop:zrank}
\emph{For all $\lm$ we have} rank$(\lm)\leq \mathrm{zrank}(\lm)$.
\end{proposition}

\textbf{Proof.} We have (see \cite[Prop.\ 7.8.3]{ec2})
  $$ h_i(1^t) = {t+i-1\choose i} =\frac{t(t+1)\cdots(t+i-1)}{i!}. $$
Hence by the Jacobi-Trudi identity,
  \beq s_{\lm}(1^t) = \det\left( {t+\lambda_i-\mu_j-i+j-1\choose 
     \lambda_i-\mu_j-i+j}\right)_{i,j=1}^n. \label{eq:slmt} \eeq
By Proposition~\ref{prop1} exactly rank$(\lm)$ rows of this matrix
have every entry equal either to 0 or a polynomial divisible by
$t$. Hence $s_{\lm}(1^t)$ is divisible by $t^{\rk(\lm)}$, so
rank$(\lm)\leq\mathrm{zrank}(\lm)$ as desired.

Alternatively, we can expand $s_{\lm}$ in terms of power sums $p_\nu$
instead of complete symmetric functions $h_\nu$. If
  \beq s_{\lm} = \sum_\nu z_\nu^{-1}\chi^{\lm}(\nu)p_\nu, 
   \label{eq:chilm} \eeq
then by the Murnaghan-Nakayama rule \cite[Cor.\ 7.17.5]{ec2}
$\chi^{\lm}(\nu)=0$ unless there exists a border strip tableau of
$\lm$ of type $\nu$. By Proposition~\ref{prop1} it follows that
$\chi^{\lm}(\nu)=0$ unless $\ell(\nu)\geq \rk(\lm)$. Since $p_\nu(1^t)
= t^{\ell(\nu)}$, it again follows that $s_{\lm}(1^t)$ is divisible by  
$t^{\rk(\lm)}$. $\ \Box$

The next result establishes that $\rk(\lm)=\mathrm{zrank}(\lm)$ in
two special cases.

\begin{theorem} \label{thm1}
(a) \emph{If $\mu=\emptyset$ (so $\lm=\lambda$) then}
$\rk(\lambda)=\mathrm{zrank}(\lambda)$.

(b) \emph{If every row of} JT$_{\lm}$ \emph{that contains a 0 also
  contains a 1, then} $\rk(\lm)=\mathrm{zrank}(\lm)$. 
\end{theorem}

\textbf{Proof.} (a) A basic formula in the theory of symmetric
functions \cite[Cor.\ 7.21.4]{ec2} asserts that
  $$ s_\lambda(1^t) = \prod_{(i,j)\in\lambda}\frac{t-i+j}{h(i,j)}, $$
where $h(i,j)=\lambda_i+\lambda'_j-i-j+1$, the hook length of
$\lambda$ at $(i,j)$. Hence 
  $$ \mathrm{zrank}(\lambda) = \#\{ i\st (i,i)\in\lambda\}
       = \rk(\lambda). $$
\indent (b) Let
  $$ y(\lm) = \left( t^{-\rk(\lm)}s_{\lm}(1^t) \right)_{t=0}. $$
By Proposition~\ref{prop:zrank} $y(\lm)$ is finite (and in fact is
just the coefficient of $t^{\rk(\lm)}$ in $s_{\lm}(1^t)$), and the
assertion that rank$(\lm)=\mathrm{zrank}(\lm)$ is equivalent to
$y(\lm)\neq 0$. Now factor out $t$ from every row not containing a 1
of the matrix on the right-hand side of equation (\ref{eq:slmt}). By
Proposition~\ref{prop1} the number of such rows is
rank$(\lm)$. Divide by $t^{\rk(\lm)}$ and set $t=0$. Denote the
resulting matrix by $R_{\lm}$, so 
  $$ y(\lm) = \left. \det R_{\lm}\right|_{t=0}. $$ 
Note that
  \beq \left(t^{-1}h_i(1^t)\right)_{t=0} = \frac 1i,\quad i\geq 1. 
   \label{eq:hit} \eeq

If row $i$ of JT$_{\lm}$ contains a 1, say in column $j$, then row
$i$ of $R_{\lm}$ has all entries equal to 0 except for a 1 in
column $j$. Hence we can remove row $i$ and column $j$ from
$R_{\lm}$ without changing the determinant $\det R_{\lm}$, except
possibly for the sign. When we do this for all rows $i$ of
JT$_{\lm}$ containing a 1, then using (\ref{eq:hit}) we obtain a
matrix of the form
  \beq R'_{\lm} = \left( \frac{1}{a_i+b_j}\right)_{i,j=1}^r, 
  \label{eq:rp} \eeq
where $a_1>a_2>\cdots>a_r>0$ and $0=b_1<b_2<\cdots<b_r$. In
particular, the denominators $a_i+b_j$ are never 0. But it was shown
by Cauchy (e.g., \cite[{\S}353]{muir}) that
  $$ \det R'_{\lm} = \frac{\prod_{i<j}(a_i-a_j)( b_i-b_j)}
        {\prod_{i,j}(a_i+b_j)}\neq 0, $$
as was to be shown. $\ \Box$

\section{Minimal border strip decompositions of $\lambda/\mu$}
\label{sec4} 
In the proof of Proposition~\ref{prop:zrank} we mentioned the
Murnaghan-Nakayama rule \cite[Cor.\ 7.17.5]{ec2} in connection with
the expansion of $s_{\lm}$ in terms of power sums. This rule asserts
that if $\chi^{\lm}(\nu)$ is defined by equation (\ref{eq:chilm}),
then 
  \beq \chi^{\lm}(\nu) = \sum_{\bms{T}} (-1)^{\mathrm{ht}(\bms{T})}, 
   \label{eq:mnrule} \eeq
summed over all border-strip tableaux $\bm{T}$ of shape $\lm$ and type
$\nu$. Here 
  $$ \mathrm{ht}(\bmt) = \sum_B \mathrm{ht}(B), $$
where $B$ ranges over all border strips in $\bm{T}$ and ht$(B)$ is one
less than the number of rows of $B$. In fact, in equation
(\ref{eq:mnrule}) $\nu$ can be composition rather than
just a partition. In other words, let $\alpha=(\alpha_1,
\dots,\alpha_m)$ be a composition of $N=|\lm|$ and let  
  $$ \chi^{\lm}(\alpha)=\sum_{\bms{T}} (-1)^{\mathrm{ht}(\bms{T})}, $$
summed over all border strip tableaux $\bm{T}$ of shape $\lm$ and type
$\alpha$. Then $\chi^{\lm}(\alpha)=\chi^{\lm}(\nu)$, where $\nu$ is
the decreasing rearrangement of $\alpha$. The second proof of
Proposition~\ref{prop:zrank} showed that $s_{\lm}$ has minimal degree
$r=\rk(\lm)$ as a polynomial in the $p_i$'s (with $\deg p_i=1$ for
$i\geq 1$). Since $p_\alpha(1^t)=t^{\ell(\alpha)}$ we see that the
coefficient $y(\lm)$ of $t^{\rk(\lm)}$ in $s_{\lm}(1^t)$ is given by
  \beq y(\lm) = \sum_{{\nu\vdash N\atop \ell(\nu)=r}} z_{\nu}^{-1}
       \chi^{\lm}(\nu). \label{eq:clm} \eeq
As mentioned above, an affirmative answer to (\ref{eq:mainop}) is
equivalent to $y(\lm)\neq 0$. Although we are unable to resolve this
question here, we will show that there is some interesting
combinatorics associated with minimal border strip decompositions and
border tableaux of shape $\lm$. In particular, a more combinatorial
version of equation (\ref{eq:clm}) is given by (\ref{eq:clmis}).

\begin{figure}
\centerline{\psfig{figure=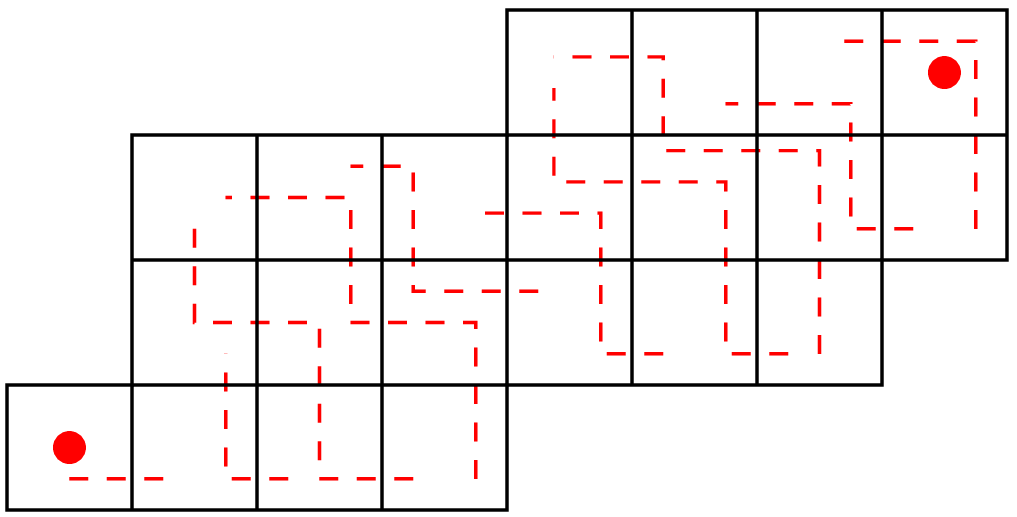}} 
\caption{Snakes for the skew shape $8874/411$}
\label{fig:snakes}
\end{figure}

Let $e$ be an edge of the lower envelope of $\lm$, i.e., no square of
$\lm$ has $e$ as its upper or left-hand edge. We will define a certain
subset $S_e$ of squares of $\lm$, called a \emph{snake}. If $e$ is
also an edge 
of the upper envelope of $\lm$, then set $S_e=\emptyset$. Otherwise,
if $e$ is horizontal and $(i,j)$ is the square of $\lm$ having $e$ as
its lower edge, then define
  \beq S_e = (\lm)\cap \{ (i,j), (i-1,j), (i-1,j-1), (i-2,j-1),
     (i-2,j-2),\dots\}. \label{eq:rsnake} \eeq
Finally if $e$ is vertical and $(i,j)$ is the square of $\lm$ having
$e$ as its right-hand edge, then define
  \beq S_e = (\lm)\cap \{ (i,j), (i,j-1), (i-1,j-1), (i-1,j-2),
     (i-2,j-2),\dots\}. \label{eq:lsnake} \eeq
In Figure~\ref{fig:snakes} the nonempty snakes of the skew shape
$8744/411$ are shown with dashed paths through their squares, with a
single bullet in the two snakes with just one square. The
\emph{length} $\ell(S)$ of a snake $S$ is one fewer than its number of
squares; a snake of length $k+1$ (so with $k$ squares) is called
a $k$-\emph{snake}. In particular, if $S_e=\emptyset$ then
$\ell(S_e)=-1$. 
Call a snake of
even length a \emph{right snake} if it has the form (\ref{eq:rsnake})
and a \emph{left snake} if it has the form (\ref{eq:lsnake}). (We
could just as well make the same definitions for snakes of odd length,
but we only need the definitions for those of even length.)  It is clear that the snakes are linearly ordered
from lower left to upper right. In this linear ordering
replace a left snake of length $2k$ with
the symbol $L_k$, a right snake of length $2k$ with $R_k$, and a snake
of odd length with $O$. The resulting sequence (which does \emph{not}
determine $\lm$), with infinitely many initial and final $O$'s
removed, is called the \emph{snake sequence} of $\lm$, denoted 
SS$(\lm)$. For instance, from Figure~\ref{fig:snakes} we see that
  $$ \sse(8874/411) = L_0OL_1L_2R_2OOL_2R_2OR_1R_0. $$

Snakes (though not with that name) appear in the solution to
\cite[Exercise~7.66]{ec2}. Call two consecutive squares of a snake $S$
(i.e., two squares with a common edge) a \emph{link} of $S$. Thus a
$k$-snake has $k-1$ links. A link of a left snake is called a
\emph{left link}, and similarly a link of a right snake is called a
\emph{right link}. Two links $l_1$ and $l_2$ are said to be
\emph{consecutive} if they have a square in common. We say that a
border strip $B$ \emph{uses} a link $l$ of some snake if $B$ contains
the two squares of $l$. Similarly a border strip decomposition $\bmd$
or border strip tableau $\bmt$ \emph{uses} $l$ if some border strip in
$\bmd$ or $\bmt$ uses $l$. The exercise cited above shows the
following. 

\begin{lemma} \label{lemma:ec2}
\emph{Let ${\cal D}$ be a border strip decomposition of $\lm$. Then no
  $B\in {\cal D}$ uses two consecutive links of a
  snake. Conversely, if we choose a set ${\cal L}$ of
  links from the snakes of $\lm$ such that no two of these links are
  consecutive, then there is a unique border strip
  decomposition ${\cal D}$ of $\lm$ that uses precisely the links in
  ${\cal L}$ (and no other links).}
\end{lemma}

Lemma~\ref{lemma:ec2} sets up a bijection between border strip
decompositions of $\lm$ and sets ${\cal L}$ of links of the snakes of $\lm$
such that no two links are consecutive. In particular, if $F_n$
denotes a Fibonacci number ($F_1=F_2=1$, $F_{n+1}=F_n+F_{n-1}$ for
$n>1$), then there are $F_{k+1}$ ways to choose a subset ${\cal L}$ of links
of a $k$-snake such that no two links are consecutive. Hence if the
snakes of $\lm$ have sizes $a_1,\dots,a_r$, then the number of border
strip decompositions of $\lm$ is $F_{a_1+1}\cdots F_{a_r+1}$ (as is
clear from the solution to \cite[Exer.\ 7.66]{ec2}). Moreover, the size
(number of border strips) of the border strip
decomposition ${\cal D}$ is given by 
 \beq \#{\cal D}=|\lm|-\#{\cal L}. \label{eq:cardd} \eeq
\indent Consider now the \emph{minimal} border strip decompositions
${\cal D}$ of $\lm$, i.e., $\#{\cal D}$ is minimized. Thus by
Proposition~\ref{prop1} we have $\#{\cal D}=\rk(\lm)$. By equation
(\ref{eq:cardd}) we wish to maximize the number of links, no two
consecutive. For snakes with an odd number $2m-1$ of links
we have no choice --- there is a unique way to choose $m$ links, no
two consecutive, and this is the maximum number possible. For 
snakes with an even number $2m$ of links there are $m+1$ ways to
choose the maximum number $m$ of links. 
Thus if mbsd$(\lm)$ denotes the number of  minimal border strip
decompositions of $\lm$, then we have proved the following result
(which will be improved in Theorem~\ref{thm:count}). 

\begin{proposition} \label{prop:nobsd}
\emph{We have}
  $$ \mathrm{mbsd}(\lm) = \prod_S \left( 1+\frac{\ell(S)}{2}\right),
  $$
\emph{where $S$ ranges over all snakes of $\lm$ of even length.}
\end{proposition}

To proceed further with the structure of the minimal border strip
decompositions of $\lm$, we will develop their connection with
code$(\lm)$. Let $p$ be the bottom-leftmost point of (the diagram of)
$\lm$, and let $q$ be the top-rightmost point. We regard the boundary
of $\lm$ as consisting of two lattice paths from $p$ to $q$ with steps
$(1,0)$ or $(0,1)$, or in other words, the restriction of the upper
and lower envelopes of $\lm$ between $p$ and $q$. 
The top-left path (regarded as a sequence of edges
$e_1,\dots,e_k$) is denoted $\Lambda_1(\lm)$, and the
bottom-right path $f_1,\dots,f_k$ by $\Lambda_2(\lm)$. Note that if in the
two-line array
  $$ \begin{array}{cccc} f_1 & f_2 & \cdots & f_k\\
                         e_1 & e_2 & \cdots & e_k \end{array} $$
we replace each vertical edge by 1 and each horizontal edge by 0, then
we obtain $\rco(\lm)$.

\begin{figure}
\centerline{\psfig{figure=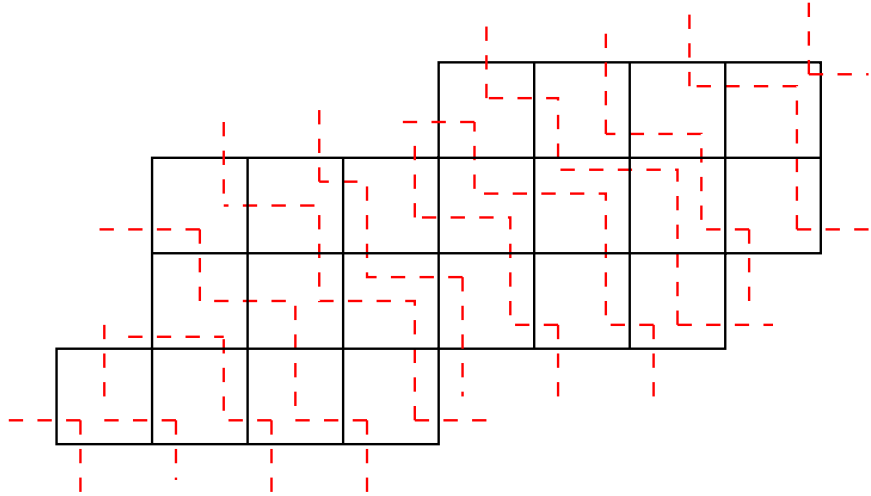}}
\caption{Extended links for the skew shape $8874/411$}
\label{fig:extend}
\end{figure}

\begin{figure}
\centerline{\psfig{figure=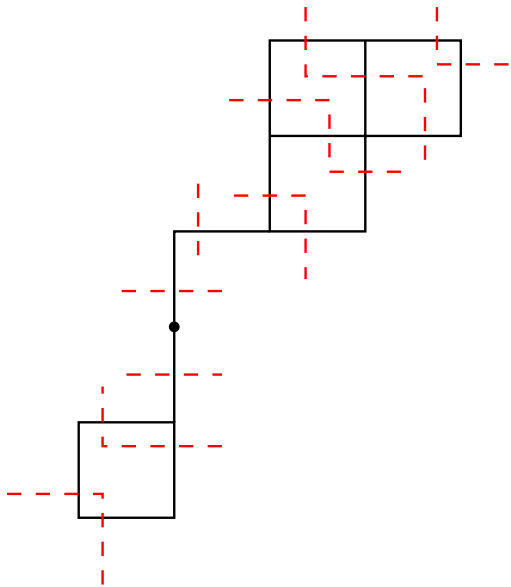}}
\caption{Extended links for the skew shape $43111/2211$}
\label{fig:dis}
\end{figure}

Continue the zigzag pattern of the links of each snake of $\lm$ one
further step in each direction, as illustrated in
Figure~\ref{fig:extend} for $\lm=8874/411$. These steps will cross an
edge on the boundary of $\lm$. Denote the top-left boundary edge
crossed by the extended link of the snake $S$ by $\tau(S)$, called the
\emph{top edge} of $S$. Similary denote the bottom-right boundary edge
crossed by the extended link of the snake $S$ by $\beta(S)$, called
the \emph{bottom edge} of the snake $S$. (In fact, the snake $S_e$ has
$\beta(S_e)=e$.) 
When $S_e=\emptyset$ we have $\tau(S_e)=\beta(S_e)=e$. 
See Figure~\ref{fig:dis} for the case
$\lm=43111/2211$, which has three edges $e$ for which $S_e=\emptyset$.

We thus have the following situation. 
Write $S_i$ as short for $S_{f_i}$, so $\tau(S_i)=e_i$ and
$\beta(S_i)=f_i$. Let
  \beq \rco(\lm) = \left( \begin{array}{cccc} c_1 & c_2 & \cdots & c_k\\
      d_1 & d_2 & \cdots & d_k \end{array} \right). \label{eq:codex}
  \eeq
It is easy to see that  $S_i$ is a left snake if and only if
$(c_i,d_i)=(1,0)$. In this case, if $S_i$ has length $2m$ then
  \beq m+1 = \#\{ j>i\st (c_j,d_j)=(0,1)\}-\#\{ j>i\st (c_j,d_j)=
  (1,0)\}. \label{eq:lls} \eeq
Similarly $S_i$ is a right snake if and only if $(c_i,d_i)=(0,1)$; and
if $S_i$ has length $2m$ then
  \beq m+1 = \#\{ j<i\st (c_j,d_j)=(1,0)\}-\#\{j<i\st (c_j,d_j)=
      (0,1)\}. \label{eq:lrs} \eeq

\begin{proposition} \label{prop:wp}
\emph{The snake sequence $\sse(\lm)=q_1q_2\cdots q_k$ is
  ``well-parenthesized'' in the following sense. There exists a
  (unique) set ${\cal P}(\lm)=\{ (u_1,v_1),\dots,(u_r,v_r)\}$, where
  $r=\rk(\lm)$, such that:}
\be
  \item[(a)] \emph{The $u_i$'s and $v_i$'s are distinct integers.}
  \item[(b)] $1\leq u_i<v_i\leq k$
  \item[(c)] \emph{$q_{u_i}=L_t$ and $q_{v_i}=R_t$ for some $t$ (depending
      on $i$)}
  \item[(d)] \emph{For no $i$ and $j$ do we have $u_i<u_j<v_i<v_j$.}
\ee
  \end{proposition}
\textbf{Proof.} Equations (\ref{eq:lp}) and (\ref{eq:lp2}) assert that
  for any $1\leq i\leq k$ we have 
  \beq \#\{ j\st 1\leq j\leq i,\ q_j=L_s\ \mbox{for some}\ s\} \geq
   \#\{ j\st 1\leq j\leq i,\ q_j=R_s\ \mbox{for some}\ s\}, 
   \label{eq:lperm} \eeq
and that the total number of $L$'s in $\sse(\lm)$ equals the total
number of $R$'s. It now follows from a standard
bijection (e.g., \cite[solution to Exer.\ 6.19(n) and (o)]{ec2}) that
there is a unique set ${\cal P}(\lm)$ satisfying (a), (b), and
(d). But (c) is then a consequence of equations (\ref{eq:lls}) and
(\ref{eq:lrs}). $\ \Box$

\begin{figure}
\centerline{\psfig{figure=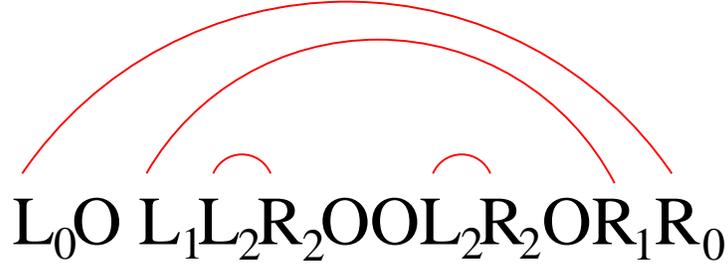}}
\caption{Parenthesization of the snake sequence $\sse(8874/411)$}
\label{fig:arcs}
\end{figure}

We can depict the set ${\cal P}(\lm)$ by drawing arcs above the terms
of $\sse(\lm)$, such that the left and right endpoints of an arc are
some $L_t$ and $R_t$, and such that the arcs are
noncrossing. For instance,
  $$ {\cal P}(8874/411) = \{(1,12),(3,11),(4,5),(8,9)\}, $$
as illustrated in Figure~\ref{fig:arcs}.

\begin{figure}
\centerline{\psfig{figure=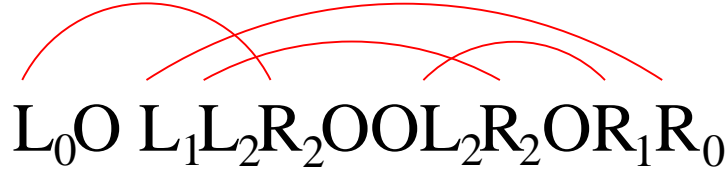}}
\caption{An interval set of the skew shape $8874/411$}
\label{fig:is}
\end{figure}

Let $\sse(\lm)=q_1q_2\cdots q_k$ as in Proposition~\ref{prop:wp}, and
define an \emph{interval set} of $\lm$ to be collection ${\cal I}$ of
$r$ ordered pairs, 
  $$ {\cal I}=\{(u_1,v_1),\dots,(u_r,v_r)\}, $$
satisfying the following conditions:
\begin{itemize}
  \item \emph{The $u_i$'s and $v_i$'s are distinct integers.}
  \item $1\leq u_i<v_i\leq k$
  \item \emph{$q_{u_i}=L_s$ and $q_{v_i}=R_t$ for some $s$ and $t$
      (depending on $i$)}
\end{itemize}
Thus ${\cal P}(\lm)$ is itself an interval set. Figure~\ref{fig:is}
illustrates the interval set $\{(1,5),(3,12),(4,9),(8,11)\}$ of the
skew shape $8874/411$. Let is$(\lm)$ denote the number of interval
sets of $\lm$. 

\begin{theorem} \label{thm:nois}
\emph{Let $T_1,\dots,T_r$ be the left snakes (or right
  snakes) of $\lm$. Then}
  $$ \mathrm{is}(\lm) = \prod_{i=1}^r \left( 1+\frac{\ell(T_i)}{2}
    \right). $$
\end{theorem}
 
\textbf{Proof.} Let $\sse(\lm)=q_1q_2\cdots q_k$. Let $q_{u_1},\dots,
q_{u_r}$ be the positions of the terms $L_s$, with
$u_1<\cdots<u_r$. Let $q_{u_i}=L_{m_i}$.
We can obtain an interval set by pairing $q_{u_r}$
with some $R_s$ to the right of $q_{u_r}$, then pairing $q_{u_{r-1}}$
with some $R_s$ to the right of $q_{u_{r-1}}$ not already paired, etc.  
By equation (\ref{eq:lls}) the number of choices for pairing $q_{u_i}$
is just $m_i+1$, and the proof follows. $\ \Box$

We are now in a position to count the number of minimal border strip
decompositions and minimal border strip tableaux of shape $\lm$. Let
us denote this latter number by mbst$(\lm)$.

\begin{theorem} \label{thm:count}
\emph{Let $\rk(\lm)=r$. Then}
  \bea \mathrm{mbsd}(\lm) & = & \mathrm{is}(\lm)^2 \label{eq:is2} \\
    \mathrm{mbst}(\lm) & = & r!\,\mathrm{is}(\lm). \label{eq:rfis} \eea
\end{theorem}

\textbf{Proof.} Equation (\ref{eq:is2}) is an immediate consequence of
Proposition~\ref{prop:nobsd} and Theorem~\ref{thm:nois} (using that in
Theorem~\ref{thm:nois} we can take $T_1,\cdots,T_r$ to consist of
either all left snakes or all right snakes).

To prove equation (\ref{eq:rfis}) we use
Proposition~\ref{prop:codebs}. Let  
   $$ \rco(\lm) = \begin{array}{cccc} 
  c_1 & c_2 & \cdots & c_k\\
  d_1 & d_2 & \cdots & d_k\end{array} $$
and let $r=\rk(\lm)$. It follows from
Proposition~\ref{prop:codebs} that a minimal border strip tableau of
shape $\lm$ is equivalent to choosing a sequence $(u_1,v_1),\dots$,
$(u_r,v_r)$ where $1\leq u_i<v_i\leq k$, $c_{u_i}=1$, $c_{v_i}=0$, the
$u_i$'s and $v_i$'s are distinct, and then successively changing
$(u_i,v_i)$ from $(1,0)$ to $(0,1)$, so that at the end we obtain the
sequence $d_1,\dots,d_k$. Since there are exactly $r$ pairs
$(c_i,d_i)$ equal to $(0,1)$ and $r$ pairs equal to $(1,0)$, the
condition that we end up with $d_1,\dots,d_k$ is equivalent to
$d_{u_i}=0$ and $d_{v_i}=1$. Hence the possible sets $\{ (u_1,v_1),
\dots, (u_r,v_r)\}$ are just the interval sets of $\lm$. There are
is$(\lm)$ ways to choose an interval set and $r!$ ways to linearly
order its elements, so the proof follows. $\ \Box$

As discussed in the above proof, every interval set ${\cal I}$ of
$\lm$ gives rise to $r!$ minimal border strip tableaux $\bmt$ of shape
$\lm$. The set of border strips appearing in such a tableau is a
border strip decomposition $\bmd$ of $\lm$. Extending our terminology
that $\bmt$ and $\bmd$ correspond to each other, we will say that
${\cal I}$, $\bmd$, and $\bmt$ all \emph{correspond} to each other.

How many of the above $r!$ border strip decompositions corresponding
to ${\cal I}$ are distinct? Rather remarkably, the number is
is$(\lm)$, independent of the interval set ${\cal I}$. This is a
consequence of Theorem~\ref{thm:isex} below. Our proof of this result
is best understood in the context of posets. Let $P$ be a finite poset
with $p$ elements $x_1,\dots,x_p$. A bijection $f:P\rightarrow
[p]=\{1,2,\dots,p\}$ is called a \emph{dropless labeling} of $P$ if we
never have $f^{-1}(i+1)<f^{-1}(i)$. Let inc$(P)$ denote the
incomparability graph of $P$, i.e, the vertex set of inc$(P)$ is
$\{x_1,\dots,x_p\}$, with an edge between $x_i$ and $x_j$ if and only
if $x_i$ and $x_j$ are incomparable in $P$. The next result is
implicit in \cite[Thm.\ 2]{g-j-w} and \cite[Theorem on p.\ 322]{b-g}
(namely, in \cite[Th.\ 2]{g-j-w} put $x=-1$ and in \cite[Theorem on
p.\ 322]{b-g} put $\lambda=-1$, and use (\ref{eq:ao}) below)
and explicit in \cite[Thm.\ 4.12]{steing}. 
For the sake of completeness we repeat the essence of the proof in
\cite{steing}. 

\begin{lemma} \label{lemma:nodrops}
  \emph{The number} dl$(P)$ \emph{of dropless labelings of $P$
    is equal to the number} ao$(P)$ \emph{of
    acyclic orientations of} inc$(P)$.
\end{lemma}

\textbf{Proof.} Given
the dropless labeling $f:P\rightarrow [p]$, define an acyclic 
orientation $\fo=\fo(f)$ as follows. If $x_ix_j$ is an edge of
inc$(P)$, then let $x_i\rightarrow x_j$ in $\fo$ if $f(x_i)<f(x_j)$,
and let $x_j\rightarrow x_i$ otherwise. Clearly $\fo$ is an acyclic
orientation of inc$(P)$. Conversely, let $\fo$ be an acyclic
orientation of inc$(P)$. The set of sources (i.e., vertices
with no arrows into them) form a chain in $P$ since otherwise two are
incomparable, so there is an arrow between them that must point into
one of them. Let $x$ be the minimal element of this chain, i.e., the
unique minimal source. If $f$ is a dropless labeling of $P$ with
$\fo=\fo(f)$, then we claim $f(x)=1$. Suppose to the contrary that
$f(x)=i>1$. Let $j$ be the largest integer satisfying $j<i$ and
$y:=f^{-1}(j)\not <x$. Note that $j$ exists since $f^{-1}(1)>x$. We
must have $y>x$ since $x$ is a source. But then $f^{-1}(j+1)\leq
x<y=f^{-1}(j)$, contradicting the fact that $f$ is dropless. Thus we
can set $f(x)=1$, remove $x$ from inc$(P)$, and proceed inductively to
construct a unique $f$ satisfying $\fo=\fo(f)$. $\ \Box$

Now given any set 
  \beq {\cal I}=\{(u_1,v_1),\dots,(u_r,v_r)\} \label{eq:intset} \eeq
with $u_i<v_i$, define a partial order $P_{{\cal I}}$ on ${\cal I}$
by setting $(u_i,v_i)<(u_j,v_j)$ if $v_i<u_j$. If we regard the
pairs $(u_i,v_i)$ as closed intervals $[u_i,v_i]$ in $\mathbb{R}$,
then $P_{{\cal I}}$ is just the \emph{interval order} corresponding
to these intervals (e.g., \cite{fish}\cite{trotter}).

\begin{lemma} \label{lemma:io}
\emph{Let ${\cal I}$ be as in equation (\ref{eq:intset}). For $1\leq
  i\leq r$ let} 
  $$ \varphi(i) = \#\{j\st v_j>v_i\}-\#\{j\st u_j>v_i\}. $$
\emph{Then} 
  $$ \mathrm{dl}(P_{{\cal I}}) = (\varphi(1)+1)(\varphi(2)+1)\cdots
       (\varphi(r)+1). $$
\end{lemma}

\textbf{Proof.} Let $\chi_{{\cal I}}(q)$ denote the chromatic
polynomial of the graph inc$(P_{{\cal I}})$. We may suppose that the
elements of ${\cal I}$ are indexed so that $v_1>v_2>\cdots>v_r$. We
can properly color the vertices of inc$(P_{{\cal I}})$ (i.e., adjacent
vertices have different colors) in $q$ colors as follows. First color
vertex $(u_1,v_1)$ in $q$ ways. Suppose that vertices $(u_1,v_1),
\dots, (u_i,v_i)$ have been colored, where $i<r$. Now for $1\leq j\leq
i$, $(u_{i+1},v_{i+1})$ is incomparable in $P_{{\cal I}}$ to
$(u_j,v_j)$ if and only $v_{i+1}>u_j$. These vertices $(u_j,v_j)$ form
an antichain in $P_{{\cal I}}$; else either some $v_j<v_{i+1}$ or some
$u_j>v_{i+1}$. The number of these vertices is $\varphi(i+1)$. Since
they form a a clique in inc$(P_{{\cal I}})$ there are exactly
$q-\varphi(i+1)$ ways to color vertex $(u_{i+1},v_{i+1})$, independent
of the colors previously assigned. It follows that
  $$ \chi_{{\cal I}}(q) = \prod_{i=1}^r (q-\varphi(i+1)). $$
For any graph $G$ with $r$ vertices it is known \cite{rs:ao} that
   \beq \mathrm{ao}(G)=(-1)^r\chi_G(-1). \label{eq:ao} \eeq
Hence
  $$ \mathrm{ao}(\mathrm{inc}(P_{{\cal I}})) =
   \prod_{i=1}^r (\varphi(i)+1). $$
The proof follows from Lemmas~\ref{lemma:nodrops} and
\ref{lemma:io}. $\ \Box$ 

\textsc{Note.} The fact (shown in the above proof) that we can order
the vertices of inc$(P_{{\cal I}})$ so that each vertex is adjacent to
a set of previous vertices forming a clique is equivalent to the
statement that the incomparability graph of an interval order is
\emph{chordal}. Note that the above proof shows that for any interval
order $P$ coming from intervals $[u_1,v_1],\dots,[u_r,v_r]$, the
chromatic polynomial of inc$(P)$ depends only on the sets $\{u_1,
\dots, u_r\}$ and $\{v_1,\dots,v_r\}$.
 
We now come to the result mentioned in the paragraph before
Lemma~\ref{lemma:nodrops}. 

\begin{theorem} \label{thm:isex}
\emph{Let ${\cal I}$ be an interval set of $\lm$, thus giving rise to
  $r!$ minimal border strip tableaux of shape $\lm$. Then the number of
  distinct border strip decompositions that correspond to these $r!$
  border strip tableaux is} is$(\lm)$.
\end{theorem}

\textbf{Proof.} Let $(u_i,v_i),(u_j,v_j)\in{\cal I}$. We say that
$(u_i,v_i)$ and $(u_j,v_j)$ \emph{overlap} if $[u_i,v_i]\cap
[u_j,v_j]\neq \emptyset$, where $[a,b]=\{u_i,u_i+1,\dots,v_i\}$.
Two linear orderings $\pi$ and $\sigma$ of ${\cal
I}$ correspond to the same border strip decomposition if and only if
any two overlapping elements $(u_i,v_i)$ and $(u_j,v_j)$ appear in the
same order in $\pi$ and $\sigma$. Suppose that $\pi$ is given by the
linear ordering 
  \beq \pi=[(u_{i_1},v_{i_1}),\dots,(u_{i_r},v_{i_r})]. \label{eq:pi} 
  \eeq 
If $(u_{i_m},v_{i_m})$ and $(u_{i_{m+1}},v_{i_{m+1}})$ are consecutive
terms of $\pi$ which do not overlap and if $i_m>i_{m+1}$, then we can
transpose the two terms without affecting the border strip
decomposition defined by $\pi$. By a series of such transpositions we
can put $\pi$ in the ``canonical form'' where consecutive
nonoverlapping pairs appear in increasing order of their
subscripts. The number of distinct border strip decompositions that
correspond to the $r!$ permutations $\pi$ is the number of $\pi$ that
are in canonical form. Let $\pi$ be given by (\ref{eq:pi}), and define
$f:P_{{\cal I}} \rightarrow [r]$ by $f(u_{i_m},v_{i_m})=m$. Then $\pi$
is in canonical form if and only if $f$ is dropless. Comparing
equation (\ref{eq:lls}), Theorem~\ref{thm:nois}, and
Lemma~\ref{lemma:io} completes the proof. $\ \Box$

Note that Theorem~\ref{thm:isex} gives a refinement of equation
(\ref{eq:is2}), since we have partitioned the is$(\lm)^2$ minimal
border strip decompositions of $\lm$ into is$(\lm)$ blocks, each of
size is$(\lm)$. 

Now let ${\cal I}=\{(u_1,v_1),\dots,(u_i,v_i)\}$ be an interval set of
$\lm$. Define the \emph{type} of ${\cal I}$ to be the partition
$\sigma$ whose parts are the integers $v_1-u_1,\dots,v_r-u_r$. Hence by
Proposition~\ref{prop:codebs} $\sigma$ is also the type of any of the
border strip decompositions corresponding to ${\cal I}$. Let
is$_\sigma(\lm)$ denote the number of interval sets of $\lm$ of type
$\sigma$, and let mbsd$_\sigma(\lm)$ denote the number of minimal
border strip decompositions of $\lm$ of type $\sigma$. The following
result is a refinement of equation (\ref{eq:is2}).

\begin{corollary} \label{cor:mbsd}
\emph{Let $N=|\lm|$. For any partition $\sigma\vdash N$, we have}
   $$ \mathrm{mbsd}_\sigma(\lm) = \mathrm{is}_\sigma(\lm)
        \mathrm{is}(\lm). $$
\end{corollary}

\textbf{Proof.} Immediate consequence of Theorem~\ref{thm:isex} and
the observation above that type$({\cal I}) = \mathrm{type}(\bm{D})$
for any interval set ${\cal I}$ and border strip decomposition
$\bm{D}$ corresponding to ${\cal I}$. $\ \Box$

We can improve the above corollary by explicitly partitioning the
minimal border strip decompositions of $\lm$ into is$(\lm)$ blocks,
each of which contains exactly mbsd$_\sigma(\lm)$ border strip
decompositions of type $\sigma$. 

\begin{theorem} \label{thm:main}
\emph{For each right snake $S$ of $\lm$ fix a set $F_S$ of $\ell(S)/2$
  links 
  of $S$, no two consecutive, and let $F=\bigcup_S F_S$. Let $Q_F$ be
  the set of all minimal border strip 
  decompositions $\bm{D}$ of $\lm$ which use the links in $Q_F$. Then
  for each $\sigma\vdash N=|\lm|$, $Q_F$ contains exactly}
  is$_\sigma(\lm)$ \emph{minimal border strip decompositions of type
  $\sigma$}. 
\end{theorem}

\begin{figure}
\centerline{\psfig{figure=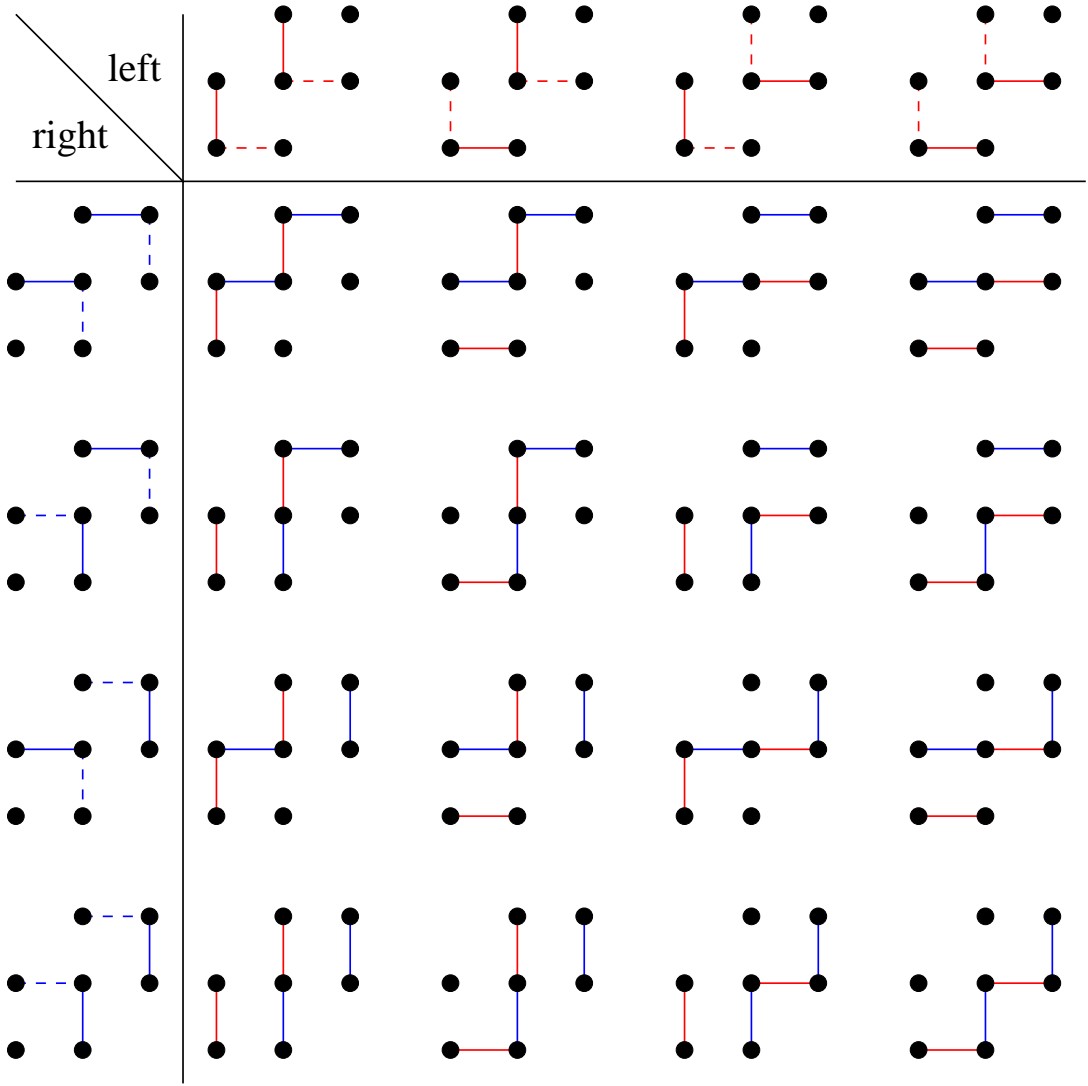}}
\caption{Minimal border strip decompositions of the skew shape $332/1$}
\label{fig:bsex}
\end{figure}

Figure~\ref{fig:bsex} illustrates Theorem~\ref{thm:main} for the case
$\lm=332/1$. We are using dots rather than squares in the diagram of
$\lm$. The first column shows the right snakes, with the choice
of links as a solid line and the remaining links as dashed lines. The
first row shows the same for the left snakes. The remaining 16 entries
are the minimal border strip decompositions of $\lm$ using the right
snake links for that row and the left snake links for that
column. Theorem~\ref{thm:main} asserts that each row (and hence by
symmetry each column) contains the same number of minimal border strip
decompositions of each type, viz., one of type $(5,1,1)$, two of type
$(4,2,1)$, and one of type $(3,2,2)$. For general $\lm$ there will
also be snakes of odd length $2m-1$ yielding $m$ links that must
be used in every minimal border strip decomposition.

\textbf{Proof of Theorem~\ref{thm:main}.} Let ${\cal I}$ be an
interval set of $\lm$ of type $\sigma$. By Theorem~\ref{thm:isex}
there are exactly is$(\lm)$ border strip decompositions (all of type
$\sigma$) corresponding to ${\cal I}$.

\emph{Claim.} Any two of the above is$(\lm)$ border strip
decompositions $\bmd$ have a different set of left links and a
different set of right links. 

\begin{figure}
\centerline{\psfig{figure=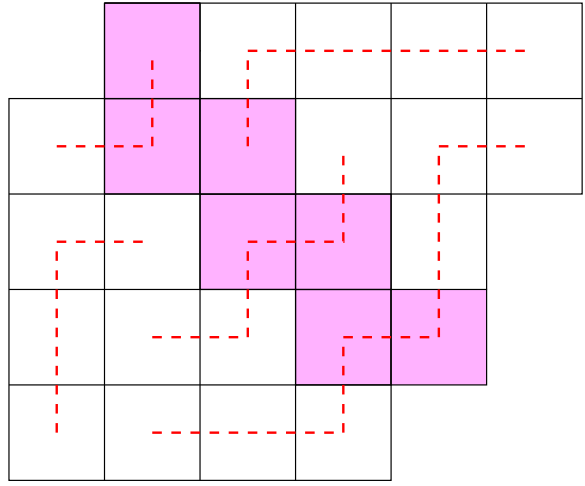}}
\caption{Intersection of border strips with a left snake}
\label{fig:intsnake}
\end{figure}

By symmetry it suffices to show that any two, say $\bmd$ and $\bmd'$,
have a different set of left links.
Let $\rco(\lm)$ be given by (\ref{eq:codex}), and let $S_i=S_{f_i}$ as
defined just before (\ref{eq:codex}). Thus $S_i$ is a left snake if
and only if $(c_i,d_i)=(0,1)$. Moreover, if $S_i$ is a left snake and
${\cal I}=\{(u_1,v_1),\dots,(u_r,v_r)\}$ is any interval set for
$\lm$, then it follows from (\ref{eq:lls}) that $\ell(S_i)=2m$ where
  $$ m = \#\{j\st u_j<i<v_j\}. $$
Let $j_1,\dots,j_m$ be those $j$ for which $u_j<i<v_j$. In a linear
ordering $\pi$ of ${\cal I}$ there are $m+1$ choices for how many of the
pairs $(u_{j_s},v_{j_s})$ precede $(u_i,v_i)$. The linear ordering
$\pi$ defines a border strip tableau with corresponding border strip
decomposition $\bmd$. In turn $\bmd$ is defined by a choice of a
maximum number of links, no two consecutive, from each left and right
snake. The choices of links from the snake $S_i$ are equivalent to
choosing the number of pairs $(u_{j_s},v_{j_s})$ preceding $(u_i,v_i)$
in $\pi$, since $S_i$ intersects precisely the border strips $B_i$ and
$B_{j_s}$ corresponding to $(u_i,v_i)$ and the $(u_{j_s},v_{j_s})$'s,
and the position of $B_i$ within the snake determines the unique two
consecutive unused links of the snake $S_i$ extended by adding one
square in each direction. Moreover, $B_i$ will be the unique
border strip whose initial square (reading from lower-left to
upper-right) begins on $S_i$. As an example see
Figure~\ref{fig:intsnake}, which shows the skew shape $\lm=66554/1$
with the left snake $S_6$ shaded. There are four border strips
intersecting $S_6$, and the third one (reading from bottom-right to
upper-left) begins on the square $(2,3)$ of $S_6$. The two links of
$S_6$ involving this square are not used in the border strip
decomposition $\bmd$.

A dropless labeling of ${\cal I}$ is uniquely determined by specifying
for each left snake $S_i$ how many of the $(u_{j_s},v_{j_s})$'s, as
defined above, precede $(u_i,v_i)$; for we can inductively determine,
preceding from left-to-right in $\rco(\lm)$, the relative order of any
pair $(u_i,v_i)$ and $(u_j,v_j)$ of elements which cross, while all
remaining ambiguities in the labeling are resolved by the dropless
condition. Thus the is$(\lm)$ dropless labelings of ${\cal I}$ define
border strip tableaux of shape $\lm$ and type $\sigma$, no two of
which have the same left links. Since these border strip tableaux
correspond to different border strip decompositions (by the proof of
Theorem~\ref{thm:isex}), the proof of the claim follows.

By the claim, for each interval set ${\cal I}$ the
is$(\lm)$ border strip decompositions corresponding to ${\cal I}$ all
have the same type and belong to different $Q_F$'s. Since there are
is$(\lm)$ different $Q_F$'s it follows that each $Q_F$ contains
exactly is$_\sigma(\lm)$ minimal border strip decompositions of type
$\sigma$, as was to be proved. $\ \Box$

Another way to state Theorem~\ref{thm:main} is as follows. Let $A$ be
the square matrix whose columns (respectively, rows) are indexed by
the maximum size sets $G$ (respectively, $F$) of links, no two
consecutive, of right snakes (respectively, left snakes) of $\lm$. The
entry $A_{FG}$ is defined to be the minimal border strip decomposition
of $\lm$ using the links $F$ and $G$. Figure~\ref{fig:bsex} shows this
matrix for $\lm=332/1$. Let $t=\mathrm{is}(\lm)$ and let ${\cal I}_1,
\dots, {\cal I}_t$ be the interval sets of $\lm$. If the border strip
decomposition $A_{FG}$ corresponds to ${\cal I}_j$, then let $L$ be
the matrix obtained by replacing $A_{FG}$ with the integer $j$. Then
the matrix $L$ is a \emph{Latin square}, i.e., every row and every
column is a permutation of $1,2,\dots,t$. For instance, when
$\lm=332/1$ the interval sets are
  $$ \begin{array}{ll}{\cal I}_1 = \{(1,6), (2,3), (4,5)\}, &
                      {\cal I}_2 = \{(1,3), (2,6), (4,5)\}\\
                      {\cal I}_3 = \{(1,5), (2,3), (4,6)\}, &
      {\cal I}_4 =\{(1,3), (2,5), (4,6)\}. \end{array} $$
The matrix $A$ of Figure~\ref{fig:bsex} becomes the Latin square
  $$ L = \left[ \begin{array}{cccc} 1 & 2 & 3 & 4\\
              2 & 1 & 4 & 3\\ 3 & 4 & 1 & 2\\ 4 & 3 & 2 & 1
          \end{array}\right]. $$

\section{An application to the characters of $\sn$.} \label{sec5}
Expand the skew Schur function $s_{\lm}$ in terms of power sums as in
equation (\ref{eq:chilm}). Define deg$(p_i)=1$, so deg$(p_\nu)=
\ell(\nu)$. As mentioned after (\ref{eq:chilm}), the Murnaghan-Nakayama
rule (\ref{eq:mnrule}) implies that if $p_\nu$ appears in $s_{\lm}$
then deg$(p_\nu)\geq r=\rk(\lm)$. In fact, at least one such $p_\nu$
actually appears in $s_{\lm}$, viz., let $\nu_1$ be the length of the
longest border strip $B_1$ of $\lm$, then $\nu_2$ the length of the
longest border strip $B_2$ of $\lm-B_1$, etc. All border strip
tableaux of $\lm$ of type $\nu$ involve the same set of border strips,
so there is no cancellation in the right-hand side of
(\ref{eq:mnrule}). Hence the coefficient of $p_\nu$ in $s_{\lm}$ in
nonzero. (See \cite[Exer.\ 7.52]{ec2} for the case $\mu=\emptyset$.)
Let us write $\hat{s}_{\lm}$ for the lowest degree part of $s_{\lm}$,
so
  \beq \hat{s}_{\lm} =\sum_{\nu\st \ell(\nu)=r}
     z_\nu^{-1}\chi^{\lm}(\nu)p_\nu, \label{eq:hats} \eeq
where $r=\rk(\lm)$. Also write  For
instance,  
  $$ s_{332/1} = \frac{1}{120}p_1^7-\frac{1}{12}p_1^4p_3
    +\frac{1}{24}p_1^3p_2^2+\frac 15 p_1^2p_5-\frac 14 p_1p_2p_4
     +\frac{1}{12}p_2^2p_3. $$
Hence
  \beas \hat{s}_{332/1} & = & \frac 15 p_1^2p_5-\frac 14 p_1p_2p_4
     +\frac{1}{12}p_2^2p_3\\[.1in] & = & \tilde{p}_1^2
      \tilde{p}_5-2\tilde{p}_1\tilde{p}_2\tilde{p}_4
        +\tilde{p}_2^2\tilde{p}_3. \eeas
If ${\cal I}=\{(w_1,y_1),\dots,(w_r,y_r)\}$ is an interval set, then
let $c({\cal I})$ denote the number of \emph{crossings} of ${\cal I}$,
i.e., the number of pairs $(i,j)$ for which
$w_i<w_j<y_i<y_j$. Moreover, let ${\cal
 P}(\lm)=\{(u_1,v_1),\dots,(u_r,v_r)\}$ be as in
Proposition~\ref{prop:wp}, and let  
  $$ \rco(\lm) = \begin{array}{cccc} 
  c_1 & c_2 & \cdots & c_k\\ 
  d_1 & d_2 & \cdots & d_k\end{array}. $$
For $1\leq i\leq r$ define
  \beas z(i) & = & \#\{ j\st u_i<j<v_i,\ c_j=0\}\\
        z(\lm) & = & z(1)+z(2)+\cdots+z(r). \eeas
It is easy to see (see the proof of Theorem~\ref{thm:topslm} for more
details) that $z(\lm)$ is just the height
$\mathrm{ht}(\bm{T})$ of a ``greedy border strip tableau'' $\bm{T}$
of shape $\lm$ obtained by starting with $\lm$ and successively
removing the largest possible border strip. (Although $\bm{T}$ may not
be unique, the set of border strips appearing in $\bm{T}$ are unique,
so ht$(\bm{T})$ is well-defined.)

\begin{lemma} \label{lemma:12}
\emph{ Let ${\cal I}$ be an interval set of $\lm$. If $\bm{T}$ and
  $\bm{T'}$ are two border strip tableaux corresponding to ${\cal I}$,
  then $\hgt(\bmt)\equiv \hgt(\bmt')\ (\mathrm{mod}\ 2)$.} 
\end{lemma}

\textbf{Proof.} When we remove a border strip $B$ of size $p$ from a
skew shape $\alpha/\beta$ with $\co(\alpha)= \cdots c_0 c_1 c_2
\cdots$, then by Proposition~\ref{prop:codebs} we replace some
$(c_i,c_{i+p})=(1,0)$ with $(0,1)$. It is easy to check (and is also
equivalent to the discussion in \cite[top of p.\ 3]{bes1}) that 
  \beq \hgt(B) = \#\{h\st i<h<i+p,\ c_h=0\}. \label{eq:htb} \eeq 
Suppose we have $(c_i,c_{i+p})=(c_j,c_{j+q})=(1,0)$, where the four
numbers $c_i,c_{i+p},c_j,c_{j+q}$ are all distinct. Let $B_1$ be the 
the border strip corresponding to $(i,i+p)$ and $B_2$ the border strip
corresponding to $(j,j+q)$ after $B_1$ has been removed. Similarly let
$B'_1$ correspond to $(j,j+q)$ and $B'_2$ to $(i,i+p)$ after $B'_1$
has been removed. If $i+p<j$ or $j+q<i$ then $B_1=B'_2$ and
$B_2=B'_1$, so $\hgt(B_1)+\hgt(B_2)=\hgt(B'_1)+\hgt(B'_2)$. In particular,
  \beq \hgt(B_1)+\hgt(B_2)\equiv\hgt(B'_1)+\hgt(B'_2)\ (\mathrm{mod}\ 2).
    \label{eq:b1b2} \eeq
If $c_i<c_j<c_{i+p}<c_{j+q}$, then using (\ref{eq:htb}) we see that 
$\hgt(B_1)=\hgt(B'_2)-1$ and $\hgt(B_2)=\hgt(B'_1)-1$ so again
(\ref{eq:b1b2}) holds. Similarly it is easy to check (\ref{eq:b1b2})
in all remaining cases. 

Iterating the above argument and using the fact that every permutation
is a product of adjacent transpositions completes the proof. $\ \Box$

\begin{theorem} \label{thm:topslm}
\emph{For any skew shape $\lm$ of rank $r$ we have}
  \beq \hat{s}_{\lm} = (-1)^{z(\lm)}\sum_{{\cal
        I}=\{(u_1,v_1),\dots,(u_r,v_r)\}} (-1)^{c({\cal I})}
       \prod_{i=1}^r \tilde{p}_{v_i-u_i}, \label{eq:top} \eeq
\emph{where ${\cal I}$ ranges over all interval sets of $\lm$.}
\end{theorem}

\textbf{Proof.} Let $\cal{I}$ be an interval set of $\lm$,
and let $\bm{T}$ be a border strip tableau corresponding to ${\cal
I}$. We claim that
  \beq \htt \equiv z(\lm)+c({\cal I})\ (\mathrm{mod}\ 2). 
    \label{eq:htt} \eeq
The proof of the claim is by induction on $c({\cal I})$.

First note that by Lemma~\ref{lemma:12}, it suffices to prove the
claim for \emph{some} $\bm{T}$ corresponding to each ${\cal I}$. 
Suppose that $c({\cal I})=0$, so ${\cal I}={\cal P}$.  Let $\bm{T}$ be
a greedy border strip tableau as defined before
Lemma~\ref{lemma:12}. The corresponding interval set is just ${\cal
P}$, the unique interval set without crossings, since if
$u_i<u_j<v_i<v_j$ we would pick the border strip corresponding to
$(u_i,v_j)$ rather than $(u_i,v_i)$ or $(u_j,v_j)$. Since by
(\ref{eq:htb}) we have $z(\lm)=\hgt(\bm{T})$, equation
(\ref{eq:htt}) holds when $c({\cal I})=0$. 

Now let $c({\cal I})>0$. Suppose that $(u_i,v_i)$ and $(u_j,v_j)$
define a crossing in ${\cal I}$, say $u_i<u_j<v_i<v_j$. Let ${\cal
  I}'$ be obtained from ${\cal I}$ by replacing $(u_i,v_i)$ and
$(u_j,v_j)$ with $(u_i,v_j)$ and $(u_j,v_i)$. It is easy to see that
$c({\cal I})-c({\cal I}')$ is an odd positive integer. By the
induction hypothesis we may assume that (\ref{eq:htt}) holds for
${\cal I}'$. Let $\bmt'$ be a border strip tableau corresponding to
${\cal I}'$ such that the border strips $B_1$ and $B_2$ indexed by
$(u_1,v_1)$ and $(u_2,v_2)$ are removed first (say in the order $B_1,
B_2$). Let $\bmt$ be the border strip tableau that differs from
$\bmt'$ by replacing $B_1,B_2$ with the border strips indexed by
$(u_j,v_i)$ and $(u_i,v_j)$. It is straightforward to verify, using
(\ref{eq:htb}) or a direct argument, that $\hgt(\bmt)$ and
$\hgt(\bmt')$ differ by an odd integer. Hence (\ref{eq:htt}) holds for
${\cal I}$, and the proof of the claim follows by induction.

Now let $\ell(\nu)=r$ and $m_i(\nu)=\#\{j\st \nu_j=i\}$, the number of
parts of $\nu$ equal to $i$. Since
$z_\nu=1^{\nu_1}\nu_1!\,2^{\nu_2}\nu_2!\cdots$, we have 
  \beas \hat{s}_{\lm} & = & \sum_{\ell(\nu)=r} z_\nu^{-1}
           \chi^{\lm}(\nu)p_\nu\\ & = & \sum_{\ell(\nu)=r}
      \frac{1}{m_1(\nu)!\,m_2(\nu)!\cdots}\chi^{\lm}(\nu)
        \tilde{p}_\nu. \eeas
Now by the Murnaghan-Nakayama rule we have
 $$ \chi^{\lm}(\nu) = \sum_{\bms{T}}(-1)^{\hgt(\bms{T})}, $$
where $\bmt$ ranges over all border strip tableaux of shape $\lm$ and
some fixed type $\alpha=(\alpha_1,\dots,\alpha_r)$ whose decreasing
rearrangement is $\nu$. Since there are $r!/m_1(\nu)!m_2(\nu)!\cdots$
different permutations $\alpha$ of the entries of $\nu$, we have
  $$ \chi^{\lm}(\nu) = \frac{m_1(\nu)!\,m_2(\nu)!\cdots}{r!}
   \sum_{\bms{T}}(-1)^{\hgt(\bms{T})}, $$
where $\bmt$ now ranges over all border strip tableaux of shape $\lm$
whose type is some permutation $\alpha$ of $\nu$. By
Theorem~\ref{thm:isex}, Proposition~\ref{prop:codebs}, and equation
(\ref{eq:htt}) we then have
  \beq \chi^{\lm}(\nu) = \frac{m_1(\nu)!\,m_2(\nu)!\cdots}{r!}
  \left( r!\sum_{{\cal I}\st \mathrm{type}({\cal I})=\nu}
   (-1)^{z(\lm)+c({\cal I})}\right), \label{eq:div} \eeq
where ${\cal I}$ ranges over all interval sets of $\lm$ of type $\nu$,
and the proof follows. $\ \Box$

Let us remark that just as in the Murnaghan-Nakayama rule,
cancellation can occur in the sum on the right-hand side of
(\ref{eq:top}). For instance, if $\lm=4442/11$ then there is one
interval set of type $(6,3,2,1)$ with one crossing and one with two
crossings.

The following corollary follows immediately from equation
(\ref{eq:div}).

\begin{corollary} \label{cor:div}
\emph{Let $\lm$ be a skew shape of rank $r$ and let
  $\ell(\nu)=r$. Then $\chi^{\lm}(\nu)$ is divisible by
  $m_1(\nu)!\,m_2(\nu)!\cdots$.}
\end{corollary}

Let $A=(a_{ij})$ be an array of real numbers with $1\leq i<j\leq
2r$. Recall that the \emph{Pfaffian} Pf$(A)$ may be defined by
(e.g.\ \cite[p.\ 616]{lov})
  $$ \mathrm{Pf}(A) =\sum_\pi (-1)^{c(\pi)} 
    a_{i_1j_1}\cdots a_{i_rj_r}, $$
where the sum is over all partitions $\pi$ of $\{1,2,\dots,2r\}$ into
two element blocks $i_k<j_k$, and where $c(\pi)$ is the number of
crossings of $\pi$, i.e., the number of pairs $h<k$ for which
$i_h<i_k<j_h<j_k$. Comparing with Theorem~\ref{thm:topslm} gives the
following alternative way of writing (\ref{eq:top}). Let $\sse(\lm)=
q_1q_2\cdots q_k$; let $u_1<u_2<\cdots<u_r$ be those indices for which
$q_{u_i}=L_s$ for some $s$; and let $v_1<v_2<\cdots<v_r$ be those
indices for which $q_{v_i}=R_s$ for some $s$. Let
$w_1<w_2<\cdots<w_{2r}$ consist of the $u_i$'s and $v_i$'s arranged in
increasing order. Then
  $$ \hat{s}_{\lm} =(-1)^{z(\lm)} \mathrm{Pf}(a_{ij}), $$
where 
  $$ a_{ij} = \left\{ \begin{array}{rl} 
   \tilde{p}_{w_j-w_i}, & \mathrm{if}\ w_i=u_s\ \mathrm{and}\ w_j=v_t\
   \mbox{for 
     some}\ s<t\\ 0, & \mbox{otherwise.} \end{array} \right. $$         
For instance, $\sse(443/2)=L_0L_1OR_1L_1R_1R_0$ and $z(443/2)=2$,
whence 
  $$ \hat{s}_{443/2} = \mathrm{Pf}\left( \begin{array}{cccccc}
     & 0 & \tilde{p}_3 & 0 & \tilde{p}_5 & \tilde{p}_6\\ & &  
    \tilde{p}_2 & 0 & \tilde{p}_4 & \tilde{p}_5\\ &  
    & & 0 & 0 & 0\\ & & & & \tilde{p}_1 & \tilde{p}_2\\ & & & & & 0\\ 
  & & & & & \end{array} \right). $$
Note that from (\ref{eq:clm}) or (\ref{eq:hats}) we get the following
Pfaffianic formula for the coefficient $y(\lm)$ of $t^{\rk(\lm)}$ in
$s_{\lm}(1^t)$:
  $$ y(\lm) = (-1)^{z(\lm)} \mathrm{Pf}(b_{ij}), $$
where 
  $$ b_{ij} = \left\{ \begin{array}{rl} 
   1/(w_j-w_i), & \mathrm{if}\ w_i=u_s\ \mathrm{and}\ w_j=v_t\
   \mbox{for 
     some}\ s<t\\ 0, & \mbox{otherwise.} \end{array} \right. $$         
Similarly from Theorem~\ref{thm:topslm} there follows
  \beq y(\lm) = (-1)^{z(\lm)}\sum_{{\cal
        I}=\{(u_1,v_1),\dots,(u_r,v_r)\}} \frac{(-1)^{c({\cal I})}}
       {\prod_{i=1}^r (v_i-u_i)}, \label{eq:clmis} \eeq
summed over all interval sets ${\cal I}$ of $\lm$.

\textsc{Acknowledgement.} I am grateful to Christine Bessenrodt for
suggesting the use of Com\'et codes in the context of skew partitions
and for supplying part (d) of Proposition~\ref{prop1}, as well as for
her careful reading of the original manuscript.

\pagebreak


\begin{thebibliography}{99}

\bibitem{bes1} C. Bessenrodt, On hooks of Young diagrams, \emph{Ann.\
Combin.}\ \textbf{2} (1998), 103--110.

\bibitem{bes2} C. Bessenrodt, On hooks of skew Young diagrams and
  bars, \emph{Ann.\ Combin.}\ \textbf{5} (2001), 37--49.

\bibitem{b-g} J. P. Buhler and R. L. Graham, A note on the binomial
  drop polynomial of a poset, \emph{J.\ Combinatorial Theory (A)}
  \textbf{66} (1994), 321--326.

\bibitem{fish} P. C. Fishburn, \emph{Interval Orders and Interval
    Graphs}, Wiley, New York, 1985.

\bibitem{g-j-w} J. R. Goldman, J. T. Joichi, and D. White, Rook theory
  III. Rook polynomials and the chromatic structure of graphs,
  \emph{J.\ Combinatorial Theory (B)} \textbf{25} (1978), 135--142.

\bibitem{lov} L. Lov\'asz, \emph{Combinatorial Problems and
    Exercises}, second ed., North-Holland, Amsterdam, 1993.

\bibitem{macd} I. G. Macdonald, \emph{Symmetric Functions and Hall
    Polynomials}, second ed., Oxford University Press, Oxford, 1995.

\bibitem{muir} T. Muir, \emph{A Treatise on the Theory of
    Determinants}, revised and enlarged by W. H. Metzler, Dover, New
    York, 1960.

\bibitem{naz} M. Nazarov and V. Tarasov, On irreducibility of tensor
  products of Yangian modules associated with skew Young diagrams,
  preprint, math.QA/0012039.

\bibitem{rs:ao} R. Stanley, Acyclic orientations of graphs,
\emph{Discrete Math.}\ \textbf{5} (1973), 171--178.

\bibitem{ec2} R. Stanley, \emph{Enumerative Combinatorics}, vol.\ 2,
  Cambridge University Press, New York/Cam\-bridge, 1999.

\bibitem{steing} E. Steingrimsson, Permutation statistics of indexed
  and poset permutations, Ph.D.\ thesis, M.I.T., 1991.

\bibitem{trotter} W. T. Trotter, \emph{Combinatorics and Partially
    Ordered Sets}, Johns Hopkins University Press, Baltimore, 1992.

\end{thebibliography}
\end{document}